\documentclass[12pt]{amsart}
\usepackage{amssymb,color,epsfig}
\usepackage[cmtip,all,ps]{xy}

\theoremstyle{plain}
\newtheorem{thm}{Theorem}[section]
\newtheorem{lem}[thm]{Lemma}

\newtheorem{prop}[thm]{Proposition}
\newtheorem{conj}[thm]{Conjecture}
\newtheorem{cor}[thm]{Corollary}

\theoremstyle{definition}
\newtheorem{dfn}[thm]{Definition}

\theoremstyle{remark}

\newtheorem{ex}[thm]{Example}
\newtheorem{rem}[thm]{Remark}

\newcommand{\ra}{\to}

\newcommand{\al}{\alpha}
\newcommand{\arrow}{\ar}

\newcommand{\bt}{\beta}

\newcommand{\Ld}{\Lambda}
\newcommand{\ga}{\gamma}
\newcommand{\fld}{\mathfrak}
\newcommand{\dt}{\mathbin{{\wr}{\wr}}}

\DeclareMathOperator{\len}{length}
\DeclareMathOperator{\Hom}{Hom}
\DeclareMathOperator{\soc}{soc}

\DeclareMathOperator{\im}{im}

\def\1#1 {\{1,\ldots,#1\} }
\def\0#1 {\{0,\ldots,#1\} }
\def\HomL{\operatorname{Hom_\Lambda}}
\def\ker{\operatorname{ker}}
\def\id{\mathrm{id}}
\def\D{\mathrm{D}}

\def\s{\mathrm{s}}
\def\t{\mathrm{t}}

\def\mod{\operatorname{mod}}
\def\Lm{\Lambda\hbox{-}\mod}
\def\Rm{R\hbox{-}\mod}
\def\disjunkt{\mathop{\hbox{$\cup$\kern - 6pt$\cdot$\kern 3pt}}}

\def\transp #1{\kern.25ex{}^t\kern-.3ex#1}
\def\NN{\mathbb{N}}

\begin{document}
\title{On Uniserial Modules in the Auslander-Reiten Quiver}
\author{Axel Boldt and Ahmad Mojiri}
\address{Department of Mathematics, Metropolitan State University,
St. Paul, MN 55106, USA}
\email{Axel.Boldt@metrostate.edu}
\address{Department of Mathematics, Texas A\&M University - Texarkana,
Texarkana, TX 75505, USA}
\email{Ahmad.Mojiri@tamut.edu}
%\address{Department of Mathematics, Wilfrid Laurier University,
%Waterloo, ON N2L 3C5, Canada}
%\email{amojiri@wlu.ca}
\thanks{The authors are grateful to Professor Walter~D.~Burgess for
his valuable suggestions and also wish to thank the referee for her/his careful reading and useful comments.}
\begin{abstract}
This article begins the study of irreducible maps involving
finite-dimensional uniserial
modules over finite-dimensional associative algebras. We work on the
classification of irreducible maps between two uniserials over
triangular algebras, and give estimates for the number of middle terms
of an almost split sequence with a uniserial end term.
\end{abstract}

\maketitle

\section{Introduction}

The study of finite-dimensional uniserial modules over
finite-di\-men\-sion\-al associative algebras was begun in earnest
by Huisgen-Zimmermann in~\cite{hz1}; Huisgen-Zimmermann and
Bongartz achieved a description of uniserial modules in terms of
certain varieties in~\cite{bhz1}. In the present article, which is based on
the authors' theses~\cite{ab} and~\cite{m}, certain
questions regarding the position of uniserial modules in the
Auslander-Reiten quiver of finite-dimensional algebras are
investigated;
most of the work applies to basic split triangular algebras only.

The article is organized as follows. In Section 2 we fix our
notation and conventions and recall the basic description of
uniserials via varieties. In Section 3 we present a general result
that motivates much of the following work: any irreducible map
between two uniserials is either the radical embedding or the
socle factor projection of a uniserial module. The two cases being
dual, we go on to state a conjecture giving a concrete necessary
and sufficient condition for a uniserial over a triangular algebra
to have an irreducible radical embedding. The criterion is
combinatorial in nature -- as a consequence, while slightly technical
when phrased in full generality, it is readily checkable for a given
quiver with relations. The sufficiency of this
condition is proved using the technique of quiver representations.
The necessity of one part of the condition is then proved in a
slightly more general context.

We have not yet managed to prove the necessity of the full
condition for all triangular algebras. In Section 4 we prove it
under an additional assumption, which includes the case of all
triangular multiserial algebras. In Section 5 we prove it for all
monomial algebras; the condition takes on a very simple form in this
situation.

In Section 6 we study a general finite-dimensional algebra and
focus on a different circle of questions: almost split
sequences with a uniserial end term. First we give a simple
general result: any short exact sequence with uniserial end terms
has a middle term which is either indecomposable or a direct sum
of two uniserials. Then we study the number of indecomposable
middle terms in an almost split sequence ending in a uniserial
module; an upper bound is given for multiserial algebras.

\section{Notation and preliminaries}

Throughout, $\fld{K}$ will be a field and $\Gamma$ will be a
finite quiver with vertex set $\Gamma_0$ and arrow set $\Gamma_1$.
We compose arrows, paths and maps from right to left: if $p:e\to
f$ and $q:f\to g$ then $qp:e\to g$. The starting point of the path
$p$ is denoted by $\s(p)$ and its end point by $\t(p)$.
$\Lambda=\fld{K}\Gamma/I$ will be a finite-dimensional
$\fld{K}$-algebra presented as the quotient of the path algebra of
$\Gamma$ by an admissible ideal $I$. $\Lambda$ is called
\emph{triangular} if $\Gamma$ does not contain any directed
cycles. Whenever useful, we identify elements of $\Gamma_0$ and
paths in $\Gamma$ with their corresponding classes in $\Lambda$.

The category of finitely generated left $\Lambda$-modules is
denoted by $\Lm$. The direct sum of two modules $M$ and $N$ is
denoted by $M \sqcup N$. A module is called \emph{uniserial} if it
has only one composition series with simple factors. If $U\in\Lm$
is uniserial with length $n$, then there exists a path $p$ in
$\Gamma$ of length $n-1$ and an element $x\in U$ such that $px\neq
0$. Any such path is called a \emph{mast} of $U$ and any such
element $x$ is called a \emph{top element} of $U$. The terminology
is that of~\cite{hz1}.

Let $p$
be a path in $\Gamma$. A path $u$ is a
\emph{right subpath} of $p$ if there exists a path $r$ with $p=ru$.
Following~\cite{hz1}, a
\emph{detour} on the path $p$ is a pair $(\al, u)$ with $\al$ an
arrow and $u$ a right subpath of $p$, where $\al u$ is a path in
$\Gamma$ which is not a right subpath of $p$, but there
exists a right subpath $v$ of $p$ with $\len(v) \geq \len(u) + 1$
such that the endpoint of $v$ coincides with the endpoint of
$\al$.
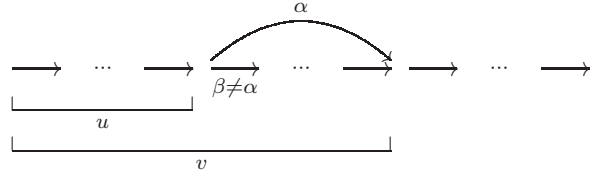
\begin{equation*}\xymatrixcolsep {1.5pc}
\begin{xy} \xymatrix {
\arrow [r] \arrow @^{|-|}[rrr]<-3ex>_{u} \arrow
@^{|-|}[0,6]<-6ex>_{v} &\arrow @{}[r]|\cdots &\arrow [r] & \arrow
[r]_{\beta\neq \alpha }\arrow @/^{1.5pc}/[rrr]^{\alpha }&\arrow
@{}[r]|\cdots &\arrow [r] & \arrow [r] &\arrow @{}[r]|\cdots
&\arrow [r] & }
\end{xy} \end{equation*}

We will abbreviate the statement ``$(\al,
u)$ is a detour on $p$'' by $(\al, u)\dt p$.
Given any detour on $p$, let
 $V(\al, u)=\{ v_i(\al, u) \mid i\in I(\al, u) \}$
be the family of right subpaths of $p$ in $\fld{K} \Gamma$ which
are longer than $u$ and have the same endpoint as $\al$.

Now suppose $p$ has length $l$ and passes consecutively through the
vertices $e(1),\dotsc ,e(l+1)$ (which need not be distinct). A \emph{route} on $p$
is any path in $\Gamma$ which starts in $e(1)$ and passes
through a subsequence of the sequence $ \left( e(1),\dotsc ,e(l+1)
\right)$ in this order and through no other vertices. A \emph{non-route} on $p$
is any path in $\Gamma$ which starts in $e(1)$ and is not a route on $p$.

Given
any uniserial module with mast $p$ and top element $x$, if $(\al,
u)\dt p$, then $\al u x =\sum_{i\in I(\al, u)} k_i(\al, u)v_i(\al,
u) x$ for  unique scalars $k_i(\al, u)$. By~\cite{hz1}, the points
$(k_i(\al, u))_{i\in I(\al, u),\, (\al, u)\dt p}$ corresponding to
uniserials with mast $p$ form an affine variety, called $V_p$,
which lives in $\mathbb{A}^N$, where $N=\sum_{(\al, u)\dt p}
|I(\al, u)|$. Moreover, there is a surjective map $\Phi_p$ from
$V_p$ onto the set of isomorphism types of uniserial $\Ld$-modules
with mast $p$. It assigns to each point $k=(k_i(\al, u))_{i\in
I(\al, u), (\al, u)\dt p}$ in $V_p$ the isomorphism type of the
module $\Ld e(1)/U_k$, where
\begin{equation*}
U_k=\sum_{(\al, u)\dt p} \Ld\left( \al u- \sum_{i\in I(\al, u)}
k_i(\al, u) v_i(\al, u)\right) +\sum_{q \text{ non-route on } p}
\Ld q.
\end{equation*}

\section{Irreducible Radical Embeddings of Uniserials}\label{s:triangembed}

In this  section, we first show that the only irreducible maps
between uniserial modules are certain radical embeddings
$JU\hookrightarrow U$ and socle factor projections $U\ra U/{\soc U}$.
Then for a triangular algebra $\Lambda=\fld{K}\Gamma/I$, we propose necessary and
sufficient combinatorial conditions for the radical embedding
$JU\hookrightarrow U$ of a uniserial module $U$ to be irreducible.

\begin{prop}\label{prop:irr} Let $R$  be a left artinian ring
with Jacobson radical $J$.
\begin{enumerate}
\item[(1)] If $f\colon M\ra U$ is an irreducible injective map from
the module $M\in\Rm$ to the uniserial $U\in\Rm$, then there exists an
isomorphism $ \varphi : JU \ra M$ so
that $f\varphi$ is the natural radical embedding $JU \hookrightarrow U$.
\item[(2)] If $g\colon U\ra M$ is an irreducible surjective map from
the uniserial $U\in\Rm$ to the module $M\in\Rm$, then there exists
an isomorphism $\psi : M \ra U/{\soc U}$
so that $ \psi g$ is the natural socle factor projection $U \ra
U/{\soc U}$.
\end{enumerate}
\end{prop}
\begin{proof}
(1) Since  $\im(f)$ is a proper submodule of $U$,
$\im(f)=J^lU$ with  $l\geq 1$ and  $M\cong J^lU$ via  $f$.
However, if $l>1$, then $J^lU\ra J^{l-1} U\ra U$ would be a
nontrivial factorization of $ J^lU\ra U$, giving us a
factorization of $f$, which is impossible. The proof of (2) is
similar to that of (1).
\end{proof}

Since every irreducible morphism is either injective or
surjective, the only irreducible maps between two uniserial modules
 are among radical embeddings $JU\hookrightarrow U$ and
socle factor projections $U\ra U/{\soc U}$. Since the two cases are
clearly dual, we will focus on radical embeddings in the sequel.

Now assume that $\Lambda=\fld{K}\Gamma/I$ is a triangular algebra. To prepare for our
analysis in this section, we fix a finitely generated uniserial
left $\Lambda$-module $U$ with mast
\[
p= 1\stackrel{\alpha_1}{\longrightarrow}
2\stackrel{\alpha_2}{\longrightarrow}\cdots\stackrel{\alpha_{n-1}}{\longrightarrow}n.
\]
On several occasions, we will refer to certain subpaths
$\alpha_i\cdots\alpha_j$ of $p$; whenever $i<j$, this expression will simply
stand for $1$. We now name all the arrows in $\Gamma$ that touch $p$,
classifying them according to the type of contact with $p$.
\begin{alignat*}{3}
&B&&:=\big\{\beta&&\in\Gamma_1\mid \text{$\s(\beta)\in\1{n-1} $ and
$\t(\beta)\not\in\1n $}\big\},\\
&B'&&:=\big\{\beta'&&\in\Gamma_1\mid \text{$\s(\beta')=n$}\big\},\\
&C&&:=\big\{\gamma&&\in\Gamma_1\mid \text{$\s(\gamma)\not\in\1{n} $ and
$\t(\gamma)\in\{2,\dots,n\}$}\big\},\\
&C'&&:=\big\{\gamma'&&\in\Gamma_1\mid \text{$\t(\gamma')=1$}\big\},\\
&D&&:=\big\{\delta&&\in\Gamma_1\mid \text{$\{\s(\delta),\t(\delta)\}\subset\1{n} $ and
$\delta\not\in\{\alpha_1,\dots,\alpha_{n-1}\}$}\big\}.
\end{alignat*}
\medbreak

For an illustration of these definitions with an example, consider the
following quiver
$\Gamma$, together with the path $p=\alpha_3\alpha_2\alpha_1$:
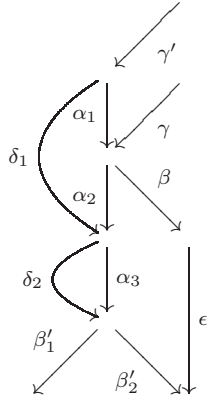
\[
\xymatrix{&&{}\ar[dl]^{\gamma'}\\
&\ar@/_5ex/[dd]_{\delta_1}\ar[d]_{\alpha_1}&{}\ar[dl]^{\gamma}\\
&\ar[d]_{\alpha_2}\ar[rd]^\beta\\
&\ar[d]^{\alpha_3}\ar@/_4ex/[d]_{\delta_2}&{}\ar[dd]^\epsilon\\
&\ar[dl]_{\beta'_1}\ar[dr]_{\beta'_2}\\
&&
}
\]
We then have
\[B=\{\beta\},\;\; B'=\{\beta_1',\beta_2'\},\;\;C=\{\gamma\},\;\;
C'=\{\gamma'\},\;\; D=\{\delta_1,\delta_2\}.\]
%\[
%\begin{array}{ll}
%B&=\{\beta\},\\
%B'&=\{\beta_1',\beta_2'\},\\
%C&=\{\gamma\},\\
%C'&=\{\gamma'\},\\
%D&=\{\delta_1,\delta_2\}.
%\end{array}
%\]

\medbreak
Observe that, in general, our uniserial module
$U$ may be identified with a representation
$U=((U_x),(f_\alpha))$ of $\Gamma$,
where
\begin{align*}
U_x&=
\begin{cases}
\fld{K},&\text{if $x\in\1n $;} \\
0,&\text{otherwise}
\end{cases}\\
\intertext{and}
f_{\alpha_i}&=\id\quad\text{for every $i\in\1{n-1} $.}
\end{align*}
The module $U$ is then completely determined by
the choice of the mast $p$ and the scalars $f_\delta(1)$ for
$\delta\in D$, different sets of scalars corresponding to
non-isomorphic modules. Unlike the hereditary case,
not every path is a mast, however,
and not every set of scalars appears in this fashion, since the
relations in $I$ impose restrictions.

We know from~\ref{prop:irr} that, in order to understand
irreducible maps between uniserial modules, it is sufficient to
study radical embeddings (and their duals, socle factor
projections). The following conjecture covers this situation; we
manage to prove ``$(2)\Rightarrow (1)$'' and a generalization of
``$(1)\Rightarrow(2)(a)$'' in the sequel. We will also prove
``$(1)\Rightarrow(2)(b)$'' for monomial and for multiserial algebras.

\begin{conj}\cite[Conjecture~1.2.1]{ab}\label{conj}
Suppose $\Ld$ is a triangular algebra and $U$ is a uniserial $\Ld$-module with mast $p$.
Then the following statements are equivalent:\\
{\rm (1)} The embedding $JU\longrightarrow U$ is irreducible.\\
{\rm (2)} $U$ is not simple and satisfies both {\rm (a)} and {\rm (b)} below:\\
\indent {\rm (a)} For every $\beta\in B $,
$$\beta\alpha_{\s(\beta)-1}\cdots\alpha_1\in Jp,$$
and for every $\delta\in D $,
$$\delta\alpha_{\s(\delta)-1}\cdots\alpha_1\in \fld{K}\alpha_{\t(\delta)-1}\cdots\alpha_1.$$
\indent {\rm (b)} There exists a subset $R\subset J$ such that
$\{rp+J^2p\mid r\in R \}$ forms a $\fld{K}$-basis for $Jp/J^2p$ and
$(i)$ and $(ii)$ both hold:\\
\indent \quad $(i)$ For every $\gamma\in C $ there exists $w\in pJ$
such that, for every $r\in R $,
$$
r\alpha_{n-1}\cdots\alpha_{\t(\gamma)}\gamma=rw.
$$
\indent \quad $(ii)$ For every $\delta\in D $ and every $r\in R $,
$$
r\alpha_{n-1}\cdots\alpha_{\t(\delta)}\delta\in
\fld{K}r\alpha_{n-1}\cdots\alpha_{\s(\delta)}.
$$
\end{conj}

\medbreak
\begin{proof}[Proof of ``$(2)\Rightarrow (1)$''.]
Let $V=((V_x),(g_\alpha))\in\Lm$ and suppose there exist
$\Lambda$-linear maps
$$
\xymatrix{JU\ar[rr]^{\Phi=(\Phi_x)}&&V\ar[rr]^{\Psi=(\Psi_x)}&&U}
$$
such that $\Psi\Phi$ is the embedding $JU\hookrightarrow U$.

Observe that we can assume without loss of generality that the
elements of the set $R$ arising from condition (2) are normed in
the following fashion: $r=e_{u(r)}re_n$ for certain vertices
$u(r)\in\Gamma_0$. We can thus denote by $g_r$ the
$\fld{K}$-linear map $V_n\longrightarrow V_{u(r)}$ induced by left
multiplication by $r$.

Note furthermore that we can strengthen the conditions on $\delta\in
D$ in the following manner:
\begin{gather*}
\delta\alpha_{\s(\delta)-1}\cdots\alpha_1=f_\delta(1)\alpha_{\t(\delta)-1}\cdots\alpha_1\\
\intertext{and for every $r\in R$}
r\alpha_{n-1}\cdots\alpha_{\t(\delta)}\delta=f_\delta(1)r\alpha_{n-1}\cdots\alpha_{\s(\delta)}.
\end{gather*}
The first equation is clear, and the second one follows then from
\[
r\alpha_{n-1}\cdots\alpha_{\t(\delta)}\delta\alpha_{\s(\delta)-1}\cdots\alpha_1=f_\delta(1)r\alpha_{n-1}\cdots\alpha_1
\]
since $rp\not=0$ for $r\in R$.\\
{\bf Case 1:} {\em There exists $v\in V_1$ with $\Psi_1(v)=1$ and
$(g_rg_{\alpha_{n-1}}\cdots g_{\alpha_1})(v)=0$ for all $r\in R $.}\\
Our goal is to construct a section $\chi$ for $\Psi$ in this case.
Define $\chi=(\chi_x):U\longrightarrow V$ by
\[
\begin{array}{rll}
\chi_1(1)&:=v,&\\
\chi_i(1)&:=(g_{\alpha_{i-1}}\cdots g_{\alpha_1})(v)&\quad\text{for
$i\in\{2,\dots,n\}$ and}\\
\chi_x&:=0&\quad\text{for $x\not\in\1n $.}
\end{array}
\]
Once we have checked that  $\chi\in\HomL(U,V)$, the equality $\Psi_1 \chi_1(1) = 1$ will
clearly imply $\Psi \chi = \id$, completing the treatment of the first case.\\
So let us check that $\chi$ is $\Lambda$-linear. That $g_{\alpha_i}\chi_i=\chi_{i+1}=\chi_{i+1}f_{\alpha_i}$ for $i\in\1{n-1} $ is
clear; moreover, we compute
\begin{align*}
g_{\delta}\chi_{\s(\delta)}(1)&=(g_{\delta}g_{\alpha_{\s(\delta)-1}}\cdots
g_{\alpha_1})(v)\\
&=f_{\delta}(1)(g_{\alpha_{\t(\delta)-1}}\cdots
g_{\alpha_1})(v)\\
&=\chi_{\t(\delta)}f_{\delta}(1)
\end{align*}
for $\delta\in D$.\\
Next observe that
$Jpv\subset \sum_r \fld{K}rpv + J^2pv
= J^2pv$ (because $R$ generates $Jp/J^2p$ and because of our
assumption in Case 1). If follows $Jpv=0$.\\
Now let $\beta\in B\cup B'$. Then
$(g_{\beta}g_{\alpha_{\s(\beta)-1}}\cdots g_{\alpha_1})(v)\in Jpv= 0$, and again
$g_{\beta} \chi_{\s(\beta)}=0=\chi_{\t(\beta)}f_\beta$.\\
\noindent{\bf Case 2:} {\em For every $v\in V_1$ with $\Psi_1(v) = 1$, there
exists $r\in R $ with \break
$(g_rg_{\alpha_{n-1}}\cdots
g_{\alpha_1})(v)\not=0$.}\\
In this case, we will construct a retraction $\chi$ for $\Phi$.
We may clearly assume that $R$ is finite, and then the condition of
Case~2 immediately implies that
there exist linear maps
$\omega_r:V_{u(r)}\longrightarrow \fld{K}$ for $r\in R $ such that
$$
\Psi_1=\sum_{r\in R}\omega_rg_rg_{\alpha_{n-1}}\cdots g_{\alpha_1}.
$$
Define $\chi=(\chi_x):V\longrightarrow JU$ by
\[
\begin{array}{rll}
\chi_i&:=\Psi_i-\sum_{r\in R}
\omega_r g_r g_{\alpha_{n-1}}\cdots g_{\alpha_i}&\quad\text{for $i\in\1n $ and}\\
\chi_x&:=0&\quad\text{for $x\not\in\1n $.}
\end{array}
\]
Again we need to check that $\chi$ is
$\Lambda$-linear.
For that purpose, we compute $\chi_1=0$,
\begin{align*}
f_{\alpha_i}\chi_i&=\Psi_{i+1}g_{\alpha_i}-(\sum_{r\in R}\omega_r g_r
g_{\alpha_{n-1}}\cdots g_{\alpha_{i+1}})g_{\alpha_i}\\
&=\chi_{i+1}g_{\alpha_i}\\
\intertext{for $i\in\1{n-1} $, and}
f_{\delta}\chi_{\s(\delta)}&=\Psi_{\t(\delta)}g_{\delta} -
f_\delta(1)\sum_{r\in R} \omega_r g_r g_{\alpha_{n-1}}\cdots
g_{\alpha_{\s(\delta)}}\\
&=\Psi_{\t(\delta)}g_{\delta} - \sum_{r\in R} \omega_r g_r g_{\alpha_{n-1}}\cdots
g_{\alpha_{\t(\delta)}}g_{\delta}\\
&=\chi_{\t(\delta)}g_{\delta}
\end{align*}
for $\delta\in D $.
In addition, we obtain
$\chi_1g_{\gamma'}=0=f_{\gamma'}\chi_{\s(\gamma')}$
for
$\gamma'\in C' $. If $\gamma\in C$, then we can clearly assume that
the corresponding element $w\in pJ$ from condition (2)(b)(i)
has the form $w=pw'$ with $w'\in
e_1Je_{\s(\gamma)}$, and it follows
\begin{align*}
\chi_{\t(\gamma)}g_{\gamma}&=f_{\gamma}\Psi_{\s(\gamma)}-\sum_{r\in
R}\omega_r g_r g_{\alpha_{n-1}}\cdots g_{\alpha_{\t(\gamma)}}g_{\gamma}\\
&=0-\sum_{r\in
R}\omega_r g_r g_{\alpha_{n-1}}\cdots g_{\alpha_1} g_{w'}\\
&=-\Psi_1g_{w'}\\
&=-f_{w'}\Psi_{\s(\gamma)}\\
&=0=f_{\gamma}\chi_{\s(\gamma)}.
\end{align*}
Hence $\chi$ belongs indeed to $\HomL(V,JU)$. That $\chi\Phi=\id_{JU}$ is a
consequence of the following computation:
\begin{align*}
\chi_2\Phi_2(1)&=\Psi_2\Phi_2(1)-\sum_{r\in R}\omega_rg_rg_{\alpha_{n-1}\cdots
\alpha_2}\Phi_2(1)\\
&=1-\sum_{r\in R} \omega_r\Phi_{u(r)} f_{r} f_{\alpha_{n-1}\cdots\alpha_2}(1)\\
&= 1\,.
\end{align*}
Thus $\Phi$ is a split monomorphism in the second case, which
shows that the inclusion $JU\hookrightarrow U$ cannot be factored
nontrivially.
\end{proof}

The implication  $(1)\Rightarrow (2)(a)$ is proved  in~\cite{ab}
using representations of algebras. In the sequel, we will
generalize $(1)\Rightarrow (2)(a)$, by weakening the assumption
that the quiver has no oriented cycle, and use the language of
modules. The following result (which does not assume that
$\Lambda$ is triangular) gives a first necessary condition for $JU
\hookrightarrow U$ to be irreducible.

\begin{prop}\label{irrembeding}
Let $U$ be a uniserial $\Lambda$-module with mast $p$. Then $JU
\hookrightarrow U$ is not irreducible if there is an arrow
leaving $e:=s(p)$ besides the first arrow of $p$.
\end{prop}
\begin{proof}
Suppose  $p=p'\bt$ with $\bt\in \Gamma_1$ and $U=\Ld e/K$
where
\[K=\sum_{(\delta,u)\dt p} \Ld\left(\delta u-\sum_{i\in I(\delta,
u)} k_i(\delta, u)v_i(\delta, u)\right)+\sum_{q\text{ non-route on
} p } \Ld q\,.\]
Suppose there is an arrow $\al$ leaving $e$ besides $\bt$.
Then either $(\al, e) \dt p$ or $\al$ is a non-route on $p$. Here, we
assume $(\al, e)\dt p$ and we will prove that $JU \hookrightarrow
U$ is not irreducible. The proof for the case where $\al$ is a
non-route on $p$ is similar.  Let $V=\Ld e/L$ with
\begin{multline*}
L= \sum_{(\delta, u)\dt p,\; (\delta, u)\neq(\al, e)}
\Ld\left(\delta u-\sum_{i\in I(\delta, u)} k_i(\delta,
u)v_i(\delta, u)\right)+\\
J\left(\alpha -\sum_{i\in I(\alpha, e)} k_i(\alpha, e) v_i(\alpha, e)\right)
 +\sum_{q\text{ non-route on } p} \Ld q\,.
\end{multline*}
We will prove that $JU\hookrightarrow U$ factors nontrivially through $V$.
Let $\varphi$ and $\psi$ be the unique $\Lambda$-homomorphisms
\[ JU \xrightarrow{\:\:\:\varphi\:\:\:} V
\xrightarrow{\:\:\:\psi\:\:\:} U\,, \] defined by
$\varphi(\bt+K)=\bt+L$ and $\psi(e+L)=e+K$. First, we show that $\varphi$ is well defined.
Suppose $\lambda \beta\in K$, for some $\lambda \in \Lambda$.
We will need to show that $\lambda \beta \in L$. We have
\[\lambda \beta =l(\alpha -\sum_{i\in I(\alpha, e)} k_i(\alpha, e) v_i(\alpha, e)) +w,\]
where $l\in\fld{K}$ and $w\in L$. On the other hand,
$\lambda \beta =k \beta +w'\beta$ for some $k\in\fld{K}$  and $w'\in J$. Hence,
\begin{equation}
l(\alpha -\sum_{i\in I(\alpha, e)} k_i(\alpha, e) v_i(\alpha, e))-k \beta=w-w'\beta\;.
\end{equation}
Since $\alpha$ does not appear
in any terms in the right hand side of (1), we have $l=0$.
Therefore $\lambda \beta=w \in L$.
It is clear that $\psi$ is well defined and $\psi\varphi$ equals the
radical embedding $JU\hookrightarrow U$.\\
Claim 1: $\varphi$ is not a split monomorphism.
Otherwise, suppose $\chi\colon V\ra JU$ is a
splitting of $\varphi$. Then
$\chi(e+L)\in JU$, thus $\chi(\bt+L)\in J^2U$, so $\bt+K=\chi\varphi(\bt +K)\in J^2 U$, which is a contradiction since $p=p'\beta$ is a mast for $U$.\\
Claim 2: $\psi$ is not a split epimorphism. Since $L\subset J\Lambda e$, the module $V$ has simple top
and is thus indecomposable. Therefore, all we have to show is that
$\psi$ is not an isomorphism. This is the case because $L$ is properly
contained in $K$ and the dimension of $V$ is larger than the dimension
of $U$.
\end{proof}

\begin{dfn}\cite{bj}
A detour $(\al, u)$ on a path $p$  is called \emph{inessential} if
\[
 \al u = s + \sum_{i\in I(\al, u)} k_i v_i (\al, u)\]
 in $\Lambda$, where $s$ is a $\fld{K}$-linear combination of paths,
 none of which is a route on $p$, and $k_i\in \fld{K}$ for all $i\in
 I(\al, u)$. A detour is \emph{essential} if it is not inessential.
\end{dfn}
The following result establishes $(1)\Rightarrow (2)(a)$ of
Conjecture~\ref{conj}, and indeed it is somewhat stronger since the
quiver is allowed to have oriented cycles here.
\begin{thm}\label{1to2a}
Let $U$ be a non-simple uniserial module with mast $p$, where $p$
does not start with an oriented cycle. If $JU \hookrightarrow U$
is irreducible, then
\begin{enumerate}
\item[(i)] All detours on $p$ are inessential.
\item[(ii)] All non-routes on $p$ are in $Jp$.
\end{enumerate}
In particular, $U=\Lambda e/{Jp}$ with $e=\s(p)$.
\end{thm}
\begin{proof}
%Recall from Proposition~\ref{irrembeding} that no detour or
%non-route starts at $e=\s(p)$.
Let $p=\al_n \cdots\al_1 $ and suppose $(\delta_j, u_j)\dt p$ for
$0\leq j \leq m$. Let $N_R=\sum_{q \text{ non-route on } p} \Lambda
q$ and
\[\Delta_j=\delta_ju_j-\sum_{i\in I(\delta_j, u_j)} k_i(\delta_j,
u_j) v_i(\delta_j, u_j)\,.\]\\
\textit{Proof of (i)}: Suppose $U=\Lambda e/K$, where $K=\sum_{j=1}^{m}
\Lambda \Delta_j + N_R$, with $m$ minimal. We have to show that
$m=0$. If $m>0$, let
$U'=\Lambda e/L$ where $L=\sum_{j=2}^{m} \Lambda \Delta_j
+J\Delta_1 + N_R$. Notice that $eJe\subseteq N_R$ by our
assumption on $p$; hence $eU'=({\fld{K} e+L})/{L}$. We first assume that $\fld{K}\neq
\mathbb{Z}_2$. Let
\[V=\frac{U'\sqcup JU'}{H}\,,\] where $H=\Lambda(p+L,
kp+L)+\Lambda(\Delta_1+L, \Delta_1+L)$ with $0,1\neq k\in \fld{K}$.
Recall from Proposition~\ref{irrembeding} that $\al_1$ is the only arrow leaving $e=\s(p)$.
Let $\varphi$ and $\psi$ be the unique $\Lambda$-homomorphisms
\[ JU \xrightarrow{\:\:\:\varphi\:\:\:} V
\xrightarrow{\:\:\:\psi\:\:\:} U\,, \] defined by
$\varphi(\al_1+K)=(\al_1+L, \al_1+L)+H$  and  $\psi((e+L, 0+L)+H)=
se+K$ and $\psi((0+L, \al_1+L)+H)=l\al_1+K$, with $s, l\in
\fld{K}$ such that $s+l=1$ and $s+lk=0$. Note that such elements
exist, since $\fld{K}\neq \mathbb{Z}_2$.
\\1. $\varphi$ is well-defined: \[\varphi(\Delta_1+K)=(\Delta_1+L, \Delta_1+L)+H=H\,.\]
2. $\psi$ is well-defined: We have $ \psi((p+L,
kp+L)+H)=sp+lkp+K=K$, and $ \psi((\Delta_1+L,
\Delta_1+L)+H)=s\Delta_1+l\Delta_1+K=K$.
\\3.  $\psi \varphi
=\id_{JU}$:
\[ \psi \varphi(\al_1 +K)=\psi((\al_1 +L,
\al_1+L)+H)=s\al_1+l\al_1+K=\al_1+K\,.\]
4. $\varphi$ is not  a split monomorphism: Otherwise there would exist $\chi \in
 \Hom_{\Lambda}(V, JU)$ such that $\chi\varphi =\id_{JU}$. Then
$\chi((e+L, L)+H)=K$, since $eJe\subseteq K$. Hence, \[\al_1+K=\chi
\varphi(\al_1+K)=\chi((\al_1+L, \al_1+L)+H)=\chi((L,
\al_1+L)+H)\,.\] Then $\chi((L, \al_1+L)+H)=\al_1+K$. Therefore,
$\chi(H)=$
\begin{eqnarray*}
\chi((p+L, kp+L)+H)&=&\chi((p+L, L)+H)+\chi((L,
kp+L)+H)\\
  &=&kp+K\neq K\,,
\end{eqnarray*}
which is a contradiction.

Therefore, $\psi$ splits; i.e.,  there
exists $\chi_1 \in \Hom_{\Ld}(U, V) $ such that $\psi \chi_1=\id_U$.
Hence $\chi_1(e+K)=(s^{-1} e+L, L)+H$ because of the assumption
that $p$ does not start with an oriented cycle. Then,
\[\chi_1(K)=\chi_1(\Delta_1+K)=\left
( s^{-1}\Delta_1+L, L \right)+H=H\, . \] Then, $(s^{-1} \Delta_1
+L, L)\in H$. Hence,
\[ (s^{-1} \Delta_1  +L, L)=z(p+L,
kp+L)+z'(\Delta_1+L, \Delta_1+L)\,,\] with $z, z'\in \Lambda$.
Therefore we have
\begin{eqnarray*}
 s^{-1} \Delta_1 +L& = & zp+z'\Delta_1 +L\,, \\
 L   & = & kzp+z'\Delta_1 +L\,.
\end{eqnarray*}
Then, $s^{-1}\Delta_1+L=(1-k)zp+L$. Hence $\Delta_1-s(1-k)zp\in
L$. Thus $s(1-k)zp\in K$, since $\Delta_1\in K$ and $L\subseteq K$.
This implies $zp\in Jp$, since $pU\neq
0$. Hence $\Delta_1 \in L$. This contradicts the minimality of $m$,
finishing the proof of (i) for these base fields.\\

Now suppose $\fld{K}=\mathbb{Z}_2$. With the same notation, let
\[V=\frac{U'\sqcup
JU'\sqcup JU'}{H}\,,\] where $H=\Lambda(L, p+L,
p+L)+\Lambda(\Delta_1+L, \Delta_1+L,\Delta_1+L)$.  Then as in the
previous case,
\[ JU \xrightarrow{\:\:\:\varphi\:\:\:} V
\xrightarrow{\:\:\:\psi\:\:\:} U\,, \] is a nontrivial
factorization of the radical embedding $JU\hookrightarrow U$ through  $V$,  where $\varphi$ and
$\psi$ are the (unique) $\Lambda$-homomorphisms defined by
$\varphi(\al_1+K)=(\al_1+L, \al_1+L,\al_1+L)+H$ and $\psi((e+L,
L, L)+H)= e+K$, $\psi((L, \al_1+L, L)+H)=\al_1+K$ and
$\psi((L, L, \al_1+L)+H)=\al_1+K$.

\textit{Proof of (ii)}:  Again first assume that $\fld{K}\neq
\mathbb{Z}_2$. By part (i), $U=\Lambda e/K$ where $K=\sum_{i=1}^m
\Lambda \bt_i u_i +Jp$, and each $\bt_i u_i$ is a non-route on $p$
with $u_i$ a right subpath of $p$, $\bt_i\in \Gamma_1$. Assume $m$
is minimum. If $m>0$, then let $U'=\Lambda e/L$ where
$L=\left({\sum_{i=2}^m \Lambda \bt_i u_i +Jp}\right)$ and
\[V=\frac{U' \sqcup JU'}{H}\,,\] where $H=\Lambda(p+L,
kp+L)+\Lambda(\bt_1 u_1+L, \bt_1 u_1+L)$ for some $k\in\fld{K},
k\neq 0,1$.
Let $\varphi$ and $\psi$ be the $\Lambda$-homomorphisms
\[ JU \xrightarrow{\:\:\:\varphi\:\:\:} V
\xrightarrow{\:\:\:\psi\:\:\:} U\,, \] defined by
$\varphi(\al_1+K)=(\al_1+L, \al_1+L)$ and  $\psi((e+L, L)+H)=s
e+K$, and $\psi((L, \al_1+L)+H)=l\al_1+K$ with $s,l\in \fld{K}$
such that  $s+l=1$ and $s+kl=0$. As in (i) we can  see that
$\varphi, \psi$ are well-defined, $\psi \varphi$ equals the radical
embedding $JU\hookrightarrow U$, and
$\varphi $ is not a split monomorphism. Therefore $\psi$ is split;
i.e., there is a $\chi\in \Hom_{\Lambda} (U, V)$ such that
$\psi\chi=\id_U$. Then $\chi(e+K)=(s^{-1} e+L, L)+H$ since no cycles
start at $p$. It follows
\[\chi(K)=\chi(\bt_1 u_1+K)= (s^{-1}\bt_1 u_1+L, L)\,.\]
Therefore,
$(s^{-1}\bt_1 u_1+L, L)=w(p+L, kp+L) + w'(\bt_1 u_1+L, \bt_1
u_1+L)$ where $w, w'\in \Lambda$. Hence,
\begin{eqnarray*}
s^{-1} \bt_1 u_1+L & = & wp+w' \bt_1 u_1+L\,,\\
  L  & = & kwp+w' \bt_1 u_1+L\,.
\end{eqnarray*}
Therefore $s^{-1}\bt_1 u_1 +L=(1-k)wp+L$. Hence
\begin{equation}\label{eqsbet}
s^{-1}\bt_1 u_1 +(k-1) wp= vp+ \sum_{i=2}^{m} w_i \bt_i u_i\,,
\end{equation}
where $v\in J$ and $w_i\in \Lambda$. If we multiply
equation~(\ref{eqsbet}) by $\t(\bt_1)$ from the left, we get that
$\t(\bt_1) wp$ is zero or a non-route on $p$, since
$\t(\bt_1)\neq \t(p)$. Then equation~(\ref{eqsbet}) contradicts the
minimality of
$m$ since it expresses $\bt_1 u_1$ as an element of $L$.\\

Now suppose that $\fld{K}=\mathbb{Z}_2$. With the same notation,
let \[V=\frac{U' \sqcup JU'\sqcup JU'}{H}\,,\] where
$H=\Lambda(L, p+L, p+L)+\Lambda(\bt_1 u_1+L, \bt_1 u_1+L, \bt_1
u_1+L)$. Let $\varphi$ and $\psi$ be the $\Lambda$-homomorphisms
\[ JU \xrightarrow{\:\:\:\varphi\:\:\:} V
\xrightarrow{\:\:\:\psi\:\:\:} U\,, \] defined by
$\varphi(\al_1+K)=(\al_1+L, \al_1+L, \al_1+L)$ and  $\psi((e+L,
L,L)+H)=e+K$,  $\psi((L, \al_1+L,   L)+H)=\al_1+K$ and $\psi((L,
L, \al_1+L)+H)=\al_1+K$. Similarly, this is a nontrivial
factorization of the radical embedding $JU\hookrightarrow U$ through $V$.
\end{proof}

\begin{ex} In order to provide a better understanding
of the different cases that would have to be dealt with in a proof of
``$(1)\Rightarrow(2)(b)$'', we include here a series of examples where
condition~(2)(b) of Conjecture~\ref{conj} is violated. A nontrivial factorization of the
radical embedding is given in each of these cases.\\
(a) Suppose $\Gamma$ is given by
$$
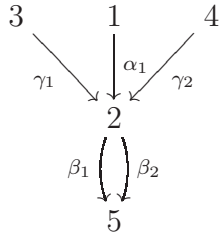
\xymatrix {3\ar[rd]_{\gamma_1}&1\ar[d]^{\alpha_1}&4\ar[ld]^{\gamma_2}\\
&2\ar@/_1ex/[d]_{\beta_1}\ar@/^1ex/[d]^{\beta_2}\\
&5}
$$
with relations $ \beta_1\alpha_1=\beta_2\alpha_1$ and
%\text{\quad and \quad}
$\beta_1\gamma_1=0=\beta_2\gamma_2.$ Here $U$ is the unique
uniserial with mast $\alpha_1$. The embedding $JU\hookrightarrow
U$ can then be factored nontrivially through a module with graph
$$
% \xymatrix{1\ar@{-}[d]&&3\ar@{-}[d]&4\ar@{-}[d]\\
% 2\ar@{-}@/_1ex/[d]\ar@{-}@/^1ex/[d]&\save .{[r].[rr]}*\frm<8pt>{-}  \restore 2\ar@{-}@/_1ex/[dl]\ar@{-}@/^1ex/[dl]&2\ar@{-}@/^1ex/[dll]^{\beta_2}&2\ar@{-}@/^3ex/[dlll]^{\beta_1}\\
% 5}
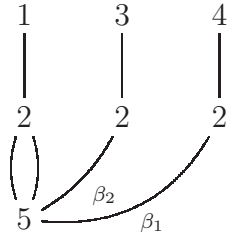
\xymatrix{1\ar@{-}[d]&3\ar@{-}[d]&4\ar@{-}[d]\\
2\ar@{-}@/_1ex/[d]\ar@{-}@/^1ex/[d]&2\ar@{-}@/^1ex/[dl]^{\beta_2}&2\ar@{-}@/^3ex/[dll]^{\beta_1}\\
5}
$$
(b) Now $\Gamma$ is given by
$$
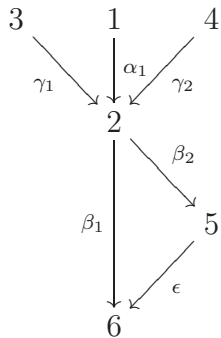
\xymatrix {
3\ar[dr]_{\gamma_1}&1\ar[d]^{\alpha_1}&4\ar[ld]^{\gamma_2}\\
&2\ar[dd]_{\beta_1}\ar[rd]^{\beta_2}\\
&&5\ar[dl]^{\epsilon}\\
&6 }$$ with relations $\epsilon\beta_2\gamma_2=\beta_1\gamma_2$
and $\beta_1\gamma_1=0=\beta_2\alpha_1.$ Again, $U$ is the unique
uniserial with mast $\alpha_1$. In this case, the radical
embedding can be factored through the indecomposable with graph
$$
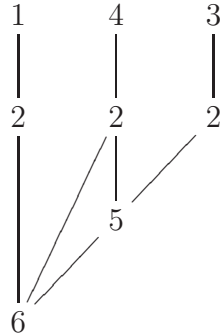
\xymatrix{
1\ar@{-}[d]&4\ar@{-}[d]&3\ar@{-}[d]\\
2\ar@{-}[dd]&2\ar@{-}[ddl]\ar@{-}[d]&2\ar@{-}[dl]\\
&5\ar@{-}[dl]\\
6
}$$
(c) Consider the quiver $\Gamma$
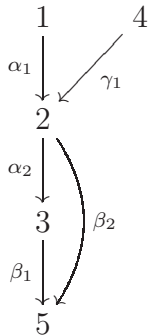
$$\xymatrix{
1\ar[d]_{\alpha_1}&4\ar[dl]^{\gamma_1}\\
2\ar[d]_{\alpha_2}\ar@/^3ex/[dd]^{\beta_2}\\
3\ar[d]_{\beta_1}\\
5 }$$ with relations $\beta_2\alpha_1=\beta_1\alpha_2\alpha_1$ and
$\beta_2\gamma_1=0.$ The radical embedding of the uniserial with
mast $\alpha_2\alpha_1$ can be factored through the following
indecomposable module:
$$
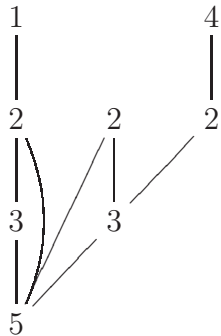
\xymatrix{
1\ar@{-}[d]&&4\ar@{-}[d]\\
2\ar@{-}@/^2ex/[dd]\ar@{-}[d]&2\ar@{-}[d]\ar@{-}[ddl]&2\ar@{-}[dl]\\
3\ar@{-}[d]&3\ar@{-}[dl]\\
5
}$$
(d) In our final example, let $\Gamma$ be given by
$$
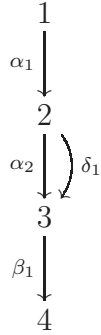
\xymatrix{
1\ar[d]_{\alpha_1}\\
2\ar[d]_{\alpha_2}\ar@/^2ex/[d]^{\delta_1}\\
3\ar[d]_{\beta_1}\\
4
}$$
and consider the relation
$
\delta_1\alpha_1=\alpha_2\alpha_1.
$
We can factor the radical embedding of the uniserial with mast
$\alpha_2\alpha_1$ through the module
$$
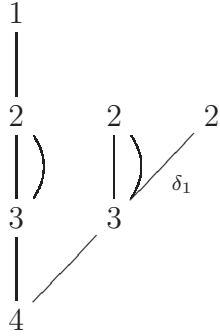
\xymatrix{
1\ar@{-}[d]\\
2\ar@{-}[d]\ar@{-}@/^2ex/[d]&2\ar@{-}[d]\ar@{-}@/^2ex/[d]&2\ar@{-}[dl]^{\delta_1}\\
3\ar@{-}[d]&3\ar@{-}[dl]\\
4
}$$
\end{ex}

\begin{rem} In order to tackle the remaining implication
``$(1)\Rightarrow(2)(b)$'' of Conjecture~\ref{conj}, it is
convenient to have the following
reformulation of condition (2)(b) at hand:\\
(2)(b') {\em There exists a family $(w_\gamma)\in(pJ)^C$, such that
for every $x\in\Gamma_0$ and $\mu\in e_xJp/e_xJ^2p$,
  we can find $r\in
e_xJe_n$ with $\mu=rp+e_xJ^2p$ and
$r\alpha_{n-1}\cdots\alpha_{\t(\gamma)}\gamma=rw_\gamma$ for all
$\gamma\in C $ and $r\alpha_{n-1}\cdots\alpha_{\t(\delta)}\delta\in
\fld{K} r\alpha_{n-1}\cdots\alpha_{\s(\delta)}$ for all $\delta\in D
$.}

Assume that condition (1) holds, i.~e., that the canonical
embedding $JU\longrightarrow U$ is irreducible, and that (2)(b') is
violated. We then get, for every family $(w_\gamma)$, a special vertex $x$ and
an element $\mu\in e_xJp/e_xJ^2p$ from the negation of this
statement. Since (2)(a)
holds, this
allows us to ``lengthen'' $U$ to a uniserial module
$\hat{U}$ in such a fashion that $U$ is an epimorphic image of
$\hat{U}$ and $\soc{\hat{U}}\simeq \Lambda e_x/J e_x$
(note however that there is a choice involved: $\hat{U}$ is
not uniquely determined by $U$ and $\mu$). Here are two potential
approaches to the construction of  a module $M$ through which the
radical embedding of $U$ factors nontrivially:\\
(a) Let $M$ be the module obtained from gluing the socles of
$\hat{U}$ and $\D(e(x)\Lambda)$ (where
$\D=\Hom_\fld{K}(-,\fld{K})$ denotes the usual duality). The
problem then is to find a ``good'' map from
$JU$ to $M$.\\
(b) This time, we begin by gluing the socles of $\hat{U}$ and
$J\hat{U}$ together to obtain $\check M$; this allows for a natural embedding of $JU$. Of course,
this particular embedding splits, and we have to extend $\check
M$ to a module $M$ having $\check M$ as an epimorphic image in order
to prevent this from happening.
\end{rem}

\section{The case of left multiserial triangular algebras}\label{s-conj-mult}

Throughout this section we assume that the algebra $\Ld$ is a
triangular algebra.
In this section, using approach (b) from above, we will show that Conjecture~\ref{conj} is true whenever the
mast $p$ has the following additional property:
\[\dim_{\fld{K}}\left( J\al_{n-1} /{J^2 \al_{n-1}}\right) \leq 1\,.\]

As a consequence, Conjecture~\ref{conj} is valid for multiserial algebras (see Corollary~\ref{mults}).
\begin{lem}\label{exist}
Let $U$ be a uniserial module with mast $p$. If the radical embedding $JU\hookrightarrow U$
is irreducible, and $\bt'$ is an arrow with $\bt' p\ne 0$, then there is a
uniserial module $V$ with mast $q:=\bt'p$. 

Specifically, if $ \{ \bt_i' p+J^2 p\; |\;
1\leq i\leq m\}$ is a basis for $Jp/{J^2p}$,
with $\bt_i'\in \Gamma_1$ and $\bt_1'=\bt'$, then such a uniserial module $V$ can be
constructed as $V=\Lambda
e/L$ with 
$$L:= Jq + \sum_{i=2}^{m} \Lambda\bt_i' p + \sum_{(\delta,
u)\dt q,\, \t(\delta)=t(q) } \Lambda \left( \delta u - l(\delta, u)
q \right)\,,$$
where $l(\delta, u)\in\fld{K}$ is a suitable scalar for
every $(\delta,u)\dt q$ with $\t(\delta)=t(q)$.
\end{lem}
\begin{proof}
%There is a basis $ \{ \bt_i' p+J^2 p\; |\; 1\leq i\leq m, \bt_i'\in
%\Gamma_1 \}$ for $Jp/{J^2p}$, with $\bt_1'=\bt'$.
Let  \[p=\: 1 \xrightarrow{\:\:\:\al_1\:\:\:} 2
\xrightarrow{\:\:\:\al_2\:\:\:} 3 \;\; \cdots
\xrightarrow{\:\:\:\al_{n-1}\:\:} n\,, \] and $n+1:=\t(\bt')$.
Suppose $(\delta, u)\dt q$. If $\t(\delta)\in \{1, 2,\dotsc, n\}$,
then $(\delta, u) \dt p$ and so by Theorem~\ref{1to2a}, $\delta u\in
\fld{K} \al_{\t(\delta)-1} \cdots \al_1 +s$, where $s\in Jp$. Since
there are no oriented cycles, we get $s=0$. If $\t(\delta)=n+1$,
then by Theorem~\ref{1to2a}(ii), $\delta u\in Jp$. Hence,
\begin{equation}\label{eqdeltau}
\delta u= l_1\bt'p+ l_2\bt_2'p+\cdots+l_m\bt_m' p + wp\,,
\end{equation}
with $w\in J^2, l_i\in \fld{K}$. Set $l(\delta, u):=l_1$. If for some $\bt \in \Gamma_1$,
$\bt u$ is a non-route on $q$, then it is a non-route on $p$ as well
and so $\bt u\in Jp$ and $\t(\bt)\notin \{1, \dotsc, n+1 \}$.
Hence, in this case, $ \bt u\in \sum_{i=2}^{m} \fld{K} \bt_i' p +
J^2 p$. Define $V=\Lambda e/L$, where
\begin{equation}\label{eqL}
L:= Jq + \sum_{i=2}^{m} \Lambda\bt_i' p + \sum_{(\delta, u)\dt
q,\, \t(\delta)=n+1 } \Lambda \left( \delta u - l(\delta,u) q \right)\,.
\end{equation}
Thus, $V$ is a uniserial module. We only need to show that $q
V\neq0$. Suppose $q V=0$.  Then, $q\in L$ and by
equations~(\ref{eqdeltau}) and~(\ref{eqL}), we get $ q\in Jq+
\sum_{i=2}^{m} \Lambda\bt_i' p + J^2 p$. Then,
\begin{equation}\label{eqex1}
q=vq + \sum_{i=2}^{m} \lambda_i \bt_i' p +  w'p\,,
\end{equation}
with $v\in J, \lambda_i\in \Lambda$ and $w'\in J^2$.  Multiply
equation~(\ref{eqex1}) by $\t(\bt')$ from the left. Since the quiver does not
have oriented cycles, $vq=0$, which  contradicts  the choice  of
the basis of $Jp/{J^2p}$.

\end{proof}
\begin{lem} \label{dim1}
Suppose $\dim_{\fld{K}} J\al_{n-1} /{J^2 \al_{n-1}}= 1$. Then
there exists an arrow $\bt'$ such that $\fld{K}\bt'\al_{n-1}
+J\bt'\al_{n-1} =J\al_{n-1}$.
\end{lem}
\begin{proof}
By the hypothesis there is some $\bt' \in \Gamma_1$ with
$\bt'\al_{n-1}\notin J^2 \al_{n-1}$. We will show that
$J^2\al_{n-1}=J\bt'\al_{n-1}$. For this we only need to show that
any path in $J^2\al_{n-1}$ is in $J\bt'\al_{n-1}$. If not, let $q$
be a longest path in $J^2\al_{n-1}\backslash J\bt'\al_{n-1}$. Then
$q=\ga_r\cdots \ga_1\al_{n-1}$, where $\ga_i \in \Gamma_1$ and $
\ga_1\al_{n-1} \notin J^2\al_{n-1}$, otherwise $q$ could be replaced
by a longer path. Hence $\ga_1 \al_{n-1}=k \bt' \al_{n-1} + w
\al_{n-1}$, where $0\neq k\in \fld{K}$ and $w\in J^2$. Therefore,
\begin{equation}\label{eqq=}
q=\ga_r\cdots \ga_1\al_{n-1}=k\ga_r\cdots
\ga_2\bt'\al_{n-1}+\ga_r\cdots \ga_2 w \al_{n-1}\,.
\end{equation}
 Since
$\ga_r\cdots \ga_2 w \al_{n-1}$ is a linear combination of paths in
$J^2 \al_{n-1}$ longer than $q$, we get $\ga_r\cdots \ga_2 w
\al_{n-1} \in J\bt'\al_{n-1}$. Then, by equation~(\ref{eqq=}), $q\in
J\bt'\al_{n-1}$. This is a contradiction.
\end{proof}
\begin{thm}\label{conj4mult}
Let $\Lambda$ be a triangular  algebra and $U$ be a uniserial
$\Lambda$-module with mast $p=\al_{n-1}\cdots\al_1$. If
$\dim_{\fld{K}} J\al_{n-1}/{J^2 \al_{n-1}}\leq 1$, then the
following statements are equivalent:
\begin{enumerate}
\item[(1)] The embedding $JU\hookrightarrow U$ is irreducible.
\item[(2)] $U$ is not simple and satisfies both (a) and (b) below:
\begin{enumerate}
\item[(a)] For every $\bt \in B$, \[ \bt \al_{\s(\bt)-1} \cdots
\al_1 \in Jp\,,\] and for every $\delta \in D$,
\[ \delta \al_{\s(\delta)-1} \cdots \al_1 \in \fld{K} \al_{\t(\delta)-1}
\cdots \al_1\,.\]
\item[(b)] $J p=0$ or there is an arrow
$\bt'$ such that $\{ \bt'p+ J^2 p \} $ forms a $\fld{K}$-basis for
$J p/{J^2 p}$ and (i) and (ii) both hold:
\begin{enumerate}
\item[(i)] For every $\gamma \in C$ there exists $w\in p J$ such
that \[ \bt' \al_{n-1} \cdots \al_{\t(\gamma)} \gamma= \bt' w\,.\]
\item[(ii)] For every $\delta \in D$ \[ \bt' \al_{n-1} \cdots
\al_{\t(\delta)} \delta \in \fld{K} \bt'\al_{n-1} \cdots
\al_{\s(\delta)}\,.
\]
\end{enumerate}
\end{enumerate}
\end{enumerate}
\end{thm}
\begin{proof}
Note first that, under the present hypotheses, the conditions (2)
are equivalent to those in Conjecture~\ref{conj}. The conditions
(2)(a) are identical. We have that $\dim_{\fld{K}} J\al_{n-1}/J^2\al_{n-1}\leq 1$ so
that, by Lemma~\ref{dim1}, we can take the set $R$ of
Conjecture~\ref{conj}(2)(b) to be $\{\bt'p+J^2p\}$ or $\emptyset$.
Then Conjecture~\ref{conj}(2)(b)(i) and (ii) reduce to the
corresponding parts of this theorem.

(1) $\Rightarrow$ (2)(b)(i):\\
Suppose $Jp\neq 0$. Then  $\dim_{\fld{K}} J\al_{n-1}/J^2\al_{n-1} =
1$. By Lemma~\ref{dim1}, there exists an arrow $\beta '\in \Gamma_1$
such that
$\fld{K}\beta'\alpha_{n-1}+J\beta'\alpha_{n-1}=J\alpha_{n-1}$. Then
$\{\bt'\al_{n-1}+J^2\al_{n-1} \}$ is a basis for
$J\al_{n-1}/{J^2\al_{n-1}}$ and $\{\bt'p+J^2p \}$ is a basis for
$Jp/{J^2p}$. We will show that for $\gamma\in C$, $ \bt' \al_{n-1}
\cdots \al_{\t(\gamma)} \gamma \in \bt' pJ$. By Theorem~\ref{1to2a},
we know that $U=\Ld e_1/{Jp}$ where $e_1=\s(p)$. Let $q=\bt'p$ and
$K= Jp$. By Lemma~\ref{exist}, there exists a uniserial module
$U_{q}=\Ld e_1/L$ with mast $q$, where
\[L=Jq+\sum_{(\delta, u)\dt q,\;\t(\delta)=\t(q)} \Lambda (\delta
u-l(\delta, u)q)\,.\]
Let
\[ V=\frac{ U_{q}\sqcup JU_{q} \sqcup \Lambda e_x}{H}\,,\] where
$H= \Lambda(q+L, q+L, 0)+ \Lambda(L, \al_{n-1} \cdots \al_1+L,
\al_{n-1} \cdots \al_{\t(\gamma)} \gamma)$ with $e_x=\s(\gamma)$.

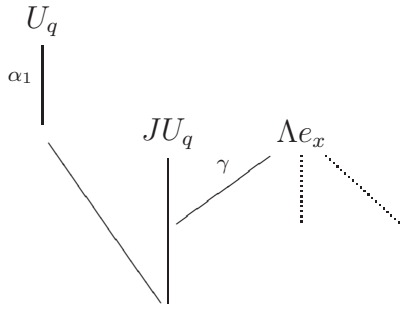
$$
\xymatrix { {U_q}\ar@{-}[d]_{\alpha_1}&& \\
\ar@{-}[ddr] & {JU_q}\ar@{-}[dd]& {\Lambda e_x}\ar@{-}[dl]_{\gamma} \ar@{.}[d] \ar@{.}[rd]\\
&&&&\\
&}
$$

%\begin{center}
%\input{pic1conj.pstex_t}
%\end{center}

Notice that for $v\in V$, $e_1v\in\fld{K}(e_1+L, 0, z)+H$, where $z$ is a linear
combination of paths from $e_x$ to $e_1$.
Let $\varphi$ and $\psi$ be the $\Lambda$-homorphisms
\[ JU \xrightarrow{\:\:\:\varphi\:\:\:} V
\xrightarrow{\:\:\:\psi\:\:\:} U\,,
\]
defined by $\varphi(\al_1+K)=(\al_1+L,\al_1+L, 0)+H$  and $\psi((e_1+L, L, 0)+H)=e_1+K,
\psi((L,\al_1+L, 0)+H)=K$ and  $\psi((L, L, e_x)+H)=K$. Then, using Lemma~\ref{dim1}
we can prove that $\varphi$ and $\psi$ are well-defined. Clearly, $\psi\varphi$
is the radical embedding $JU\hookrightarrow U$.\\
Claim:  $\varphi$ is not a split monomorphism; otherwise there would
exist $\chi\colon V\ra JU$ such that $\chi \varphi=\id_{JU}$. We have
$\al_1+K =\chi\varphi(\al_1+K)=\chi((\al_1+L, \al_1+L,
 0)+H)=\chi((\al_1+L, L, 0)+H)+\chi((L, \al_1+L, 0)+H)=
\chi((L, \al_1+L, 0)+H)$, because $\chi((e_1+L, L, 0)+H)=K$. Also we
have $\chi((L, L, e_x)+H)=K$. But
\begin{multline*}
\chi(H) =\chi((L, \al_{n-1} \cdots
\al_1+L, \al_{n-1} \cdots \al_{\t(\gamma)}\gamma)+H)\\=
\al_{n-1}\cdots\al_2\al_1+K\neq K\,,
\end{multline*}
which is a contradiction.
Therefore $\psi$ splits, i.e., there exists $\chi_1:U\ra V$ with
$\psi \chi_1=\id_U.$ We have $\chi_1(e_1+K)=((e_1+L, L,
\sum_{i=1}^{m} k_i w_i)+H),$  where $w_i$ are the paths from $e_x$
to $e_1$ and $k_i\in \fld{K}$. But $q\in K$ and so
\[\chi_1(K)=\chi_1(q+K)= (q+L, L, \sum_{i=1}^{m} k_i q w_i)+H\,.\]
Hence, $(q+L, L, \sum_{i=1}^{m}
k_i q w_i) \in $\[ \Lambda(q+L, q+L, 0) + \Lambda(L,
\al_{n-1}\cdots\al_1+L, \al_{n-1} \cdots\al_{\t(\gamma)}\gamma)\,.
\] Then, by Lemma~\ref{dim1}, $(q+L, L, \sum_{i=1}^{m} k_i q w_i)=$
\begin{multline*}
k(q+L, q+L, 0)
 + l \bt'(L, \al_{n-1}\cdots\al_1+L, \al_{n-1}
\cdots\al_{\t(\gamma)}\gamma))\\
+\sum_{l(u_i)\geq 1} l_i u_i \bt'(L, \al_{n-1}\cdots\al_1+L,
\al_{n-1} \cdots\al_{\t(\gamma)}\gamma)\,,
\end{multline*}
where $k, l, l_i\in \fld{K}$. It follows $k=1 \text{ and } l=-1$.
Hence, \[\bt'\al_{n-1} \cdots \al_{\t(\gamma)} \gamma = - \sum k_i q
w_i+\sum_{l(u_i)\geq 1} l_iu_i\bt'\al_{n-1}\cdots\al_{\t(\gamma)}
\gamma\,.\]  If we multiply the above equation from the left by
$\t(\bt')$, using the fact that the quiver does not have any
oriented cycles and therefore $\t(\bt') u_i=0$, we obtain
  \[\bt'\al_{n-1}\cdots\al_{\t(\gamma)} \gamma=-
\sum k_i q w_i \in \bt' pJ\,.\]

(1) $\Rightarrow$ (2)(b)(ii):\\
Suppose $\delta \in D$. We will show that
$\bt'\al_{n-1}\cdots\al_{\t(\delta)}\delta \in
\fld{K}\bt'\al_{n-1}\cdots\al_{\s(\delta)}$. Let $\delta \colon
i\ra j$ and
 $q:=\bt'\al_{n-1}\cdots \al_1$.  Again let
$U_{q}=\Ld e/L$ be the uniserial with mast $q$, with $L$ as above.
 Let \[ V=\frac{ U_{q}\sqcup JU_{q} \sqcup
\Lambda e_i}{H}\,,\] where $H=\Lambda(q+L, q+L, 0)+$ \[ \Lambda(L,
\al_{n-1} \cdots \al_1+L, \al_{n-1} \cdots \al_j \delta)+
\Lambda(L, L, \al_{n-1}\cdots\al_i)\,.\] Let $\varphi$ and $\psi$
be the $\Lambda$-homomorphisms
\[ JU \xrightarrow{\:\:\:\varphi\:\:\:} V
\xrightarrow{\:\:\:\psi\:\:\:} U\,,
\]
defined by $\varphi(\al_1+K)=(\al_1+L, \al_1+L, 0)+H$ and $\psi((e_1+L,
L, 0)+H)=e_1+K, \psi((L, \al_1+L, 0)+H)=K$ and  $\psi((L, L,
e_i)+H)=K$. Then, $\varphi$ and $\psi$ are well-defined and
$\psi\varphi$ is the radical embedding $JU\hookrightarrow U$.\\
Claim: $\varphi$ is not split monomorphism; otherwise there would
exist $\chi\colon V\ra JU$ such that $\chi \varphi=\id_{JU}$. Then we
would have $\al_1+K =\chi\varphi(\al_1+K )=\chi((\al_1+L, \al_1+L,
 0)+H)=\chi((\al_1+L, L, 0)+H)+\chi((L, \al_1+L, 0)+H)=
\chi((L, \al_1+L, 0)+H)$, because $\chi((e_1+L, L, 0)+H)=K$.
By Theorem~\ref{1to2a}(ii), $\chi((L, L, e_i)+H)=k\al_{i-1}\cdots \al_1+K$, where
$k\in \fld{K}$ . Thus, $\chi(H)=\chi((L, L, \al_{n-1}\cdots \al_i
)+H)=k\al_{n-1}\cdots\al_1 +K$. Therefore, $k=0$. But $\chi(H)=$
\[\chi((L, \al_{n-1}\cdots \al_1+L,
\al_{n-1}\cdots\al_j\delta)+H)= \al_{n-1}\cdots\al_2\al_1+K \neq
K\,,\] which is a contradiction.
%3- $\psi$ is well-defined.\\

Therefore, $\psi$  splits, i.e., there exists $\chi_1:U\ra V$ with
$\psi \chi_1=\id_U$. We have $\chi_1(e_1+K)=(e_1+L, L, 0)+H$. Hence
$\chi_1(K)=\chi_1(q+K)=(q_1+L, L, 0)+H.$ Therefore,
\begin{multline*}
(q+L, L, 0)\in \Lambda(q+L, q+L, 0)+ \Lambda(L,
\al_{n-1}\cdots\al_1+L, \al_{n-1}\cdots\al_j\delta)\\ +\Lambda(L,
L, \al_{n-1} \cdots\al_i)\,.
\end{multline*}
Then by Lemma~\ref{dim1}
\begin{multline*}
(q+L, L, 0) =k(q+L, q+L, 0) +l \bt'(L, \al_{n-1}\cdots\al_1+L,
 \al_{n-1}\cdots \al_j\delta)\\
 + \sum_{l(u_s)\geq1} l_s
u_s\bt'(L, \al_{n-1}\cdots\al_1+L, \al_{n-1}\cdots\al_j\delta)\\
+ v(L, L,   \al_{n-1}\cdots\al_i)\,,
\end{multline*}
 where $l, l_s \in
\fld{K}$,  $ u_s\in J$ and $ v\in \Lambda$. Hence $k=1$ and
$l=-1$. Therefore,
\begin{equation}\label{E:beta}
\bt'\al_{n-1}\cdots\al_j\delta
= \sum_{l(u_s)\geq 1} l_s u_s\bt'\al_{n-1}\cdots\al_j\delta + v
\al_{n-1}\cdots\al_i\,,
\end{equation}
in $\Lambda e_i$. If we multiply equation~(\ref{E:beta})
from the left  by $\t(\bt')$; using the fact that there are no
oriented cycles, $\t(\bt') u_i=0$.  We get
\[\bt'\al_{n-1}\cdots\al_j\delta=\t(\bt') v \al_{n-1}\cdots\al_i\,.\]
Since there are no oriented cycles, $t(\beta')\ne t(\alpha_{n-1})$
and so $\t(\bt') v\al_{n-1} \in J\al_{n-1}$. But
$J\al_{n-1}=\fld{K}\bt'\al_{n-1} + J\bt'\al_{n-1}$ by
Lemma~\ref{dim1}. Therefore,
\[\bt'\al_{n-1}\cdots\al_j\delta=k\bt'\al_{n-1}\cdots\al_i+
w\bt'\al_{n-1}\cdots\al_i\,,\] where $w\in J$. But, $\t(\bt')w=0$,
since there are no  oriented cycles. Therefore,
$\bt'\al_{n-1}\cdots\al_j\delta=k\bt'\al_{n-1}\cdots\al_i$.
\end{proof}
By the work above, Conjecture~\ref{conj} is true for all triangular
algebras with a presentation so that for each $\al\in \Gamma_1$,
$\Ld \al$ is uniserial.
\begin{dfn}
An algebra $\Ld$ with Jacobson radical $J$ is called \emph{left multiserial}
(\emph{$m$-multiserial}) if, for each primitive idempotent $e$ of
$\Ld$, the left ideal $Je$ is a sum of uniserial ($m$ uniserial)
$\Ld$-modules.
\end{dfn}
For the convenience of the reader, we provide here the following theorem from~\cite[Remark 2.3]{bj}.
\begin{thm}~\cite[Remark 2.3]{bj}
Every left multiserial algebra is isomorphic to
one with a presentation so that for each $\al\in \Gamma_1$,
$\Ld \al$ is uniserial.
\end{thm}
\begin{cor}~\label{mults}
Conjecture~\ref{conj} is true for all left triangular multiserial
algebras.
\end{cor}

\section{The case of monomial algebras}

Throughout this section we assume
that the algebra $\Ld$ is a triangular algebra. We will
prove that the conjecture is true for monomial algebras.
Recall that for any path $p$, nonzero in $\Lambda$, there is an
affine variety $V_p$ and a map $\Phi_p$ from $V_p$ onto the set of isomorphism types
of uniserial $\Lambda$-modules with mast $p$ (see page 3).

\begin{thm}
Suppose $\Lambda$ is a triangular monomial algebra  and $U$ is a
uniserial $\Lambda$-module with mast $p$. Then the following
statements are equivalent:
\begin{enumerate}
\item[(1)] The embedding $JU\hookrightarrow U$ is irreducible.
\item[(2)] $U$ is not simple and satisfies both (a) and (b) below:
\begin{enumerate}
\item[(a)]
\begin{enumerate}
\item[(i)] For every $\bt\in B$, $\bt\al_{\s(\bt)-1}\cdots\al_1=0$, and
\item[(ii)] For every $\delta \in D$, $\delta
\al_{\s(\delta)-1}\cdots\al_1=0$.
\end{enumerate}
\item[(b)] For every $\bt'\in B'$ such that $\bt'p\neq0$ we have:
\begin{enumerate}
\item[(i)] For every $\gamma \in C,$ $ \bt' \al_{n-1} \cdots
\al_{\t(\gamma)} \gamma= 0$, and
\item[(ii)] For every $\delta \in D,$
$ \bt' \al_{n-1} \cdots \al_{\t(\delta)} \delta =0$.
\end{enumerate}
\end{enumerate}
\end{enumerate}
\end{thm}
\begin{proof}
Note first that, since the algebra is monomial, the conditions (2)
are equivalent to the ones in Conjecture~\ref{conj}.

(1) $\Rightarrow$ (2)(b)(i):\\
Let $p=\al_{n-1}\cdots\al_1$ and $U=\Ld e_1/K$. Suppose that there
is $\bt'\in B'$ such that $\bt'p\neq 0$ and $\bt'
\al_{n-1}\cdots\al_i\gamma\neq0$ for some $\gamma \in C$, where $ x
\xrightarrow{\:\:\gamma\:\:} i, \text{ with } x\notin \{1,2,\dotsc
,n\}$. By condition (2)(a), $V_p=\{\underline{0}\}$. Let
\[ q_1:=\bt'\al_{n-1}\cdots\al_1,
\;q_2:=\bt'\al_{n-1}\cdots\al_i\gamma\,.\] Since $\Lambda$ is a
monomial algebra and $q_i\neq 0$;  by~\cite[Proposition II.3]{m1},
$\underline{0}\in V_{q_i}$ for  $i=1$ and  $2$. Let
$U_{q_1}:=\Phi_{q_1}(\underline{0})=\Ld e_1/L$  and
$U_{q_2}:=\Phi_{q_2}(\underline{0})=\Ld e_x/F$, where
$e_x=\s(\gamma)$.
%where $\Varphi_{q_i}\colon
%V_{q_i}\ra \{ \text{ Uniserials with mast } q_i \} $.
 Let \[ V=\frac{U_{q_1}\sqcup
JU_{q_1}\sqcup U_{q_2}}{H}\,,\] where
\[H=\Lambda(q_1+L, q_1+L, F)+\Lambda(L,
\al_{n-1}\cdots\al_1+L, \al_{n-1}\cdots\al_i\gamma+F)\,.\] Once
again, for $v\in V$, $e_1v=(ke_1+L, L, z+F)+H$, where $z$ is a
linear combination of paths from $\s(\gamma)$ to $e_1$. However,
such a path goes through $e_1$ and so is a non-route on $q_2$,
i.e., $z\in F$.  Let $\varphi$ and $\psi$ be the $\Lambda$-homomorphisms
\[ JU \xrightarrow{\:\:\:\varphi\:\:\:} V
\xrightarrow{\:\:\:\psi\:\:\:} U\,,
\]
defined by $\varphi(\al_1+K)=(\al_1+L,
\al_1+L, F)+H$ and $ \psi((e_1+L, L, F)+H)
= e_1+K,\, \psi((L, \al_1+L, F)+H)=K$ and    $\psi((L, L,
e_x+F)+H)=K$.
We will first show that $\varphi$ is well-defined. Note that $K=L+\Lambda q_1$.
Suppose $\lambda\alpha_1\in K$, for some $\lambda\in \Lambda$. Then
$\lambda\alpha_1=w+\gamma q_1$, for some $w\in L$ and $\gamma \in \Lambda$.
Thus, $(\lambda\alpha_1+L, \lambda\alpha_1+L, F)=(\gamma q_1+L, \gamma q_1+L, F)\in H$.\\
Again, $\psi$ is well-defined, $\psi \varphi$ is $JU\hookrightarrow
U$ and $\varphi$ is not split monomorphism.  We will prove that
$\psi$ also is not a split epimorphism, which contradicts the
irreducibility of $JU\hookrightarrow U$. Suppose $\psi$ is a split
epimorphism. Then, there exists $\chi:U\ra V$ with  $\psi \chi=\id_U$.
We have $\chi(e_1+K)=(e_1+L, L, F)+H$. But
$q_1=\bt'\al_{n-1}\cdots\al_1 \in K$. Hence
$\chi(K)=\chi(q_1+K)=(q_1+L, L, F)+H$ is zero in $V$. Then,
\begin{equation}\label{E:q}
\begin{split}
(q_1+L, L, F)=k&(q_1+L, q_1+L ,F)\\+&w(L, \al_{n-1}\cdots\al_1+L,
\al_{n-1}\cdots\al_{\t(\gamma)}\gamma+F)\,,
\end{split}
\end{equation}
where $k \in\fld{K}$ and $w\in \Lambda$. Note that if $\beta_1'\neq
\beta'$, then either $\beta_1'\alpha_{n-1}\cdots\alpha_1$ is
non-route on $q_1$ or $(\beta_1', \alpha_{n-1}\cdots\alpha_1)\dt
q_1$. Hence, $\beta_1'\alpha_{n-1}\cdots\alpha_1\in L$. Similarly,
$\beta_1'\alpha_{n-1}\cdots\alpha_{t(\gamma)}\gamma\in F$. Therefore
equation~(\ref{E:q}) becomes $(q_1+L, L, F) =k(q_1+L, q_1+L
,F)+l\bt'(L, \al_{n-1}\cdots\al_1+L,
\al_{n-1}\cdots\al_{\t(\gamma)}\gamma+F)\,,$ where $k,l \in
\fld{K}$. Therefore $k=1,\; k+l=0, \;l=0$, which is a contradiction.

(1) $\Rightarrow$ (2)(b)(ii):\\
Suppose there is $\bt'\in B'$ such that $\bt'p\neq 0$ and
$\bt'\al_{n-1}\cdots\al_{\t(\delta)}\delta\neq0$ for some $\delta\in
D$. Let $\delta\colon i\ra j$. By (a)(ii), $\s(\delta)=i \neq 1$.
 Let \[ q_1:=\bt'\al_{n-1}\cdots\al_1,
q_2:=\bt'\al_{n-1}\cdots\al_j\delta\,.\] and let $U_{q_1}=\Ld e_1/L$
and $U_{q_2}=\Ld e_2/F$ be the uniserial modules corresponding to
$0\in V_{q_1}$ and $0\in V_{q_2}$ respectively.  Let
\[
V=\frac{U_{q_1}\sqcup JU_{q_1}\sqcup U_{q_2}}{H}\,,\] where
$H=\Lambda(q_1+L, q_1+L, F)+\Lambda(L, \al_{n-1}\cdots\al_1+L,
\al_{n-1}\cdots\al_j\delta+F)$.

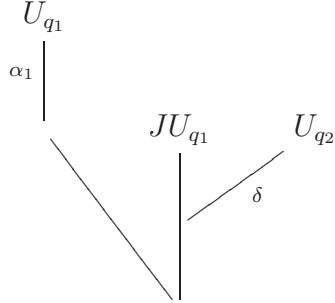
$$
\xymatrix { {U_{q_1}}\ar@{-}[d]_{\alpha_1}&& \\
\ar@{-}[ddr] & {JU_{q_1}}\ar@{-}[dd]& {U_{q_2}}\ar@{-}[dl]^{\delta}\\
&\\
&}
$$

%\begin{center}
%\input{pic2conj.pstex_t}
%\end{center}

Let $\varphi$ and $\psi$ be the $\Lambda$-homomorphisms
\[ JU \xrightarrow{\:\:\:\varphi\:\:\:} V
\xrightarrow{\:\:\:\psi\:\:\:} U\,,
\]
defined by $\varphi(\al_1+K)=(\al_1+L,
\al_1+L, F)+H$  and  $\psi((e_1+L, L,
F)+H)=e_1+K,\, \psi((L, \al_1+L, F)+H)=K$ and  $\psi((L, L,
e_i+F)+H)=K$. Again, $\varphi$
and $\psi$ are well-defined and $\psi\varphi$ is
the radical embedding $JU\hookrightarrow U$.\\
Claim: $\varphi$ is not a split monomorphism:\\
Suppose there exists $\chi\colon V\ra JU$ such that $\chi
\varphi=\id_{JU}$. Then, we have $\al_1+K =\chi\psi(\al_1+K
)=\chi((\al_1+L, \al_1+L,
 F)+H)=\chi((\al_1+L, L, F)+H)+\chi((L, \al_1+L, F)+H)=
\chi((L, \al_1+L, F)+H)$. Also we know that $\chi((L, L,
e_i+F)=k\al_{i-1}\cdots \al_1+K$, where $k\in \fld{K}$. Then
$\chi((L, L, \delta e_i+F)= k\delta\al_{i-1}\cdots \al_1+K=K$, by
(a)(ii),  and
\[\chi(H)=\chi(L, \al_{n-1}\cdots \al_1+L,
\al_{n-1}\cdots\al_j\delta +F)= \al_{n-1}\cdots\al_2\al_1 +K\neq
K\,,\] which is a
contradiction.\\
Claim: $\psi $ is not a  split epimorphism:\\
Suppose there exists $\chi_1:U\ra V$ with  $\psi \chi_1=\id_U$. We have
$\chi_1(e_1+K)=(e_1+L, L, F)+H$. Hence
$\chi_1(K)=\chi_1(q_1+K)=(q_1+L, L, F)+H$. Therefore $(q_1+L,L,
F)+H=H$, and so \[ (q_1+L, L, F)\in \Lambda(q_1+L, q_1+L, F)
+\Lambda \bt'(L, \al_{n-1}\cdots\al_1+L,
\al_{n-1}\cdots\al_j\delta+F)\,.\] Then $(q_1+L, L,
F)=k(q_1+L,q_1+L, F)+l\bt'(L, \al_{n-1}\cdots\al_1+L,
\al_{n-1}\cdots\al_j\delta+F)$, with $k,l \in \fld{K}$. Therefore
$k=1,\; k+l=0, \;l=0$, which is a contradiction.

\end{proof}

\section{Almost split sequences with uniserial end terms}

In this section, we first show that if we have an arbitrary exact
sequence with uniserial end terms, then the middle term is either
indecomposable or a direct sum of two uniserials. Then we study
$\al(U)$, the number of indecomposable summands of the middle term
of an almost split sequence ending in $U$, where $U$ is a uniserial
nonprojective $\Ld$-module and give a global upper bound for it in
the case that $\Ld$ is a multiserial algebra.
\begin{prop}
Let $R$ be a left artinian ring and consider a short exact
sequence
$$
\xymatrix{0\ar[r]&U_1\ar[r]^f&M\ar[r]^g&U_2\ar[r]&0}
$$
in $\Rm$ with uniserial modules $U_1$ and $U_2$. Then $M$ is
either indecomposable or a direct sum of two uniserial modules.
\end{prop}

\begin{proof}
We will again denote the Jacobson radical of $R$ by $J$. Assume we
have a decomposition $M=M_1\oplus M_2$ with both $M_1$ and $M_2$
non-zero. Decompose $f$ and $g$ accordingly, i.~e., write
$f={f_1\choose f_2}$ and $g=(g_1,g_2)$, and let
$$
\bar{{}}:\Rm\longrightarrow (R/J)\hbox{-}\operatorname{mod}
$$
be the functor $R/J\underset{R}{\otimes}-$. We then get the right
exact sequence
$$
\xymatrix{{\bar{U_1}}\ar[r]^-{\bar{f_1}\choose\bar{f_2}}&{\bar{M_1}}\oplus\bar{M_2}\ar[r]^-{(\bar{g_1},\bar{g_2})}&{\bar{U_2}}\ar[r]&0}
$$
where $\bar{U_1}$ and $\bar{U_2}$ are simple and $\bar{M_1}$,
$\bar{M_2}$ non-zero semisimple. Comparing the lengths of the
involved modules, we see that both $\bar{M_1}$ and $\bar{M_2}$
must be simple and $\bar{f}\not=0$. Without loss of generality, we
may assume $\bar{f_1}(\bar{U_1})=\bar{M_1}$.

Pick $u_1\in U_1\setminus JU_1$. Then $f_1(u_1)\in M_1\setminus
JM_1$ generates $M_1$. Hence $f_1$ is surjective and $M_1$ is
uniserial. If $f_2(u_1)=0$, then $f_2=0$ and $g_2$ is injective,
and consequently $M_2$ is uniserial. If $f_2(u_1)\not=0$, we can
find $l\geq 0$ with $f_2(u_1)\in J^lM_2\setminus J^{l+1}M_2$. If
$l=0$, then $f_2(u_1)$ generates $M_2$ and $M_2$ is therefore
uniserial. We will assume $l>0$ from now on.

\noindent{\bf Claim 1:} {\em $\im(g_1)\subset J^l U_2$.}

\noindent Let $m_1\in M_1$; write $m_1=\alpha f_1(u_1)=f_1(\alpha
u_1)$ with $\alpha\in\Lambda$. Then $g_1(m_1)=g(m_1)=gf_1(\alpha
u_1)-gf(\alpha u_1)=-gf_2(\alpha u_1)\subset g(J^l M_2)\subset J^l
U_2$. Hence we have $g_1(M_1)\subset J^lU_2$.

\noindent{\bf Claim 2:} {\em $g_2$ is surjective and the map
$M_2/J^lM_2\longrightarrow U_2/J^l U_2$ induced by $g_2$ is an
isomorphism.}

\noindent Let $m_2\in M_2\setminus JM_2$. Then $u_2:=g_2(m_2)\in
U_2\setminus JU_2$ (since $g_2(m_2)\in JU_2$ would imply
$\im(g)=\im(g_1)+\im(g_2)\subset J^lU_2 + JU_2\subsetneqq U_2$, a
contradiction). Since $u_2$ generates $U_2$, $g_2$ is surjective.
 Now let $x\in M_2\setminus J^l M_2$ and assume
$g_2(x)\in J^lU_2$, say $g_2(x)=\alpha u_2=g_2(\alpha m_2)$ with
$\alpha\in J^l$. Then $x-\alpha m_2\in \ker(g_2)\setminus
J^lM_2\subset \im(f_2)\setminus J^l M_2=\emptyset$, again a
contradiction.

\noindent{\bf Claim 3:} {\em $J^lM_2$ is uniserial.}

\noindent By restricting our maps $f$ and $g$, we obtain the
following short exact sequence:
$$
\xymatrix{0\ar[r]&U_1\ar[r]&M_1\oplus J^lM_2\ar[r]&J^lU_2\ar[r]&0}
$$
and we see as above that $J^lM_2/J^{l+1}M_2$ is simple, hence
$J^lM_2$ is generated by $f_2(u_1)$ and $f_2:U_1\longrightarrow
J^lM_2$ is therefore surjective.

\noindent{\bf Claim 4:} {\em $M_2$ is uniserial.}

\noindent We know that $J^kM_2/J^{k+1}M_2$ is simple or $0$ for
all $k\in\NN$.
\end{proof}

In the sequel, $\Lambda$ will be a finite-dimensional algebra over
$\fld{K}$.

The following proposition gives a general upper bound for the number
$\al(U)$ for a uniserial module $U$:

\begin{prop} \label{soc2}
If $U\in\Lm$ is a non-projective uniserial module, then \[\al(U) \leq
\len(\soc DTrU) +1\;.\]
\end{prop}
\begin{proof}
Let $0 \ra DTrU \ra B \ra U \ra 0$ be an almost split sequence.
Then  $0 \ra \soc DTrU \ra \soc B \ra \soc U $ is left exact. Therefore,
\begin{eqnarray*}
\al (U) &\leq& \len(\soc B)\\
&\leq& \len(\soc DTrU) + \len(\soc U)\\
&=& \len(\soc DTrU) +1\,.
\end{eqnarray*}
\end{proof}

The following proposition gives more precise information.

\begin{prop}\label{mono}
Let $0\ra DTrU \xrightarrow{\:\:f\:\:} \bigsqcup_{i\in
I}B_i\xrightarrow{\:\:g\:\:} U\ra 0$ be an almost split sequence
where $U$ is a uniserial module and the $B_i$ are indecomposable.
\begin{enumerate}
\item [(i)]  At most one of the induced  maps
$g_i\colon B_i\ra U$ is a monomorphism. \item [(ii)] If
$B_i\xrightarrow{\:\:g_i\:\:} U$  is an  epimorphism and $\soc
B_i$ is simple then $\soc B_i \subseteq f(\soc DTrU)$.
\item [(iii)] Let $I'=\{i\in I \mid g_i\colon
B_i\ra U \text{ is an  epimorphism}\}$\,. Then  $|I'| \leq
\len(\soc DTrU)$.
\end{enumerate}
\end{prop}
\begin{proof}
(i) Suppose $g_1$ and  $g_2$ are monomorphisms. Using
Proposition~\ref{prop:irr} again, we have $B_1 \cong
JU$ and  $B_2 \cong JU$. The induced irreducible morphism
$B_1\sqcup B_2\ra U $ cannot be an epimorphism  and therefore is a
monomorphism and $B_1\sqcup B_2\cong JU$, which is impossible. \\
   (ii) We have $\soc B_i \cap \ker(g_i)\neq 0$ since
$\ker(g_i)\neq 0$ and  $\soc B_i$ is essential in $B_i$. But $\soc
B_i$ is simple, so $\soc B_i \subseteq \ker(g_i)$. We know that
$0\ra \soc DTrU \xrightarrow{\:\:\bar{f}\:\:} \bigsqcup_{i\in
I}\soc B_i\xrightarrow{\:\:\bar{g}\:\:} \soc U$ is exact. Hence
$\soc B_i \subseteq \ker \bar{g}=\im \bar{f}$.  Therefore, $\soc
B_i \subseteq f(\soc DTrU)$.\\
    (iii) We distinguish two cases:\\
Case 1: There is an $i$ such that $g_i$ is a monomorphism. Then
\[ |I'|\leq\al (U) -1\leq \len (\soc DTrU) \]
by Proposition~\ref{soc2}.\\
Case 2: For each $i\in I$, the map $g_i$ is an epimorphism. We
consider the exact sequence $0\ra \soc DTrU
\xrightarrow{\:\:\bar{f}\:\:} \bigsqcup_{i\in I}\soc
B_i\xrightarrow{\:\:\bar{g}\:\:} \soc U$ and we use (ii): if $\soc
B_i $ is simple for all $i$, then $\bar{f}$ is an isomorphism and
we get
\[
|I'| = \al (U) = \len (\soc \bigsqcup_{i\in I} B_i) =
\len (\soc DTrU)\,.
\]
If however at least one $\soc B_i $ is not simple, then the same exact
sequence gives
\[
|I'| = \al (U) \leq \len (\soc \bigsqcup_{i\in I} B_i) -1 \leq
\len (\soc DTrU)\,.
\]
\end{proof}

Let $e, f$ be primitive idempotents in $\Lambda$. For a non-zero
element $a \in fJe$, the $\Lambda$-module $\Lambda e/ {\Lambda a}$
is indecomposable and non-projective. We are interested in
the case where this module is a uniserial module and consider the almost split
sequence ending in $\Lambda e/ {\Lambda a}$.

\begin{prop} \label{indec}
If $U=\Lambda e/{\Lambda a}$ is a uniserial module, then $\al(U)\leq 2$.
\end{prop}
\begin{proof}
$\Lambda f \xrightarrow{\:\:.a\:\:} \Lambda e\ra \Lambda e/
{\Lambda a}\ra 0$ (where $.a$ denotes the right
multiplication by $a$) is exact and is the start of a minimal projective presentation of
$\Lambda e/ {\Lambda a}$. From~\cite[Proposition V.6.1]{ars} we have
that the middle term $B$ in the almost split sequence $\delta
\colon 0\ra DTrU \ra B\ra U\ra 0$ has a decomposition $B=B'\sqcup
B''$ with $B'$ indecomposable and such that if $B''\neq 0$, the
induced morphism $g''\colon B''\ra U$ is an irreducible
monomorphism. But, by Proposition~\ref{prop:irr}, $B'' \cong JU$
is indecomposable and therefore $\al(U) \leq 2$.
\end{proof}

Uniserial representations of left multiserial algebras are studied
in \cite{bj}.  Here we find an upper bound for $\al(U)$ where $U$
is a uniserial module over a left $m$-multiserial algebra.

\begin{thm}\label{alpha2}
Let $U$ be a non-projective uniserial module over a left
$m$-multiserial algebra $\Lambda$ with $m\geq 2$. Then $\al(U)\leq
m$.
\end{thm}
\begin{proof}
By \cite[Remark 2.3]{bj}, we can assume that
$\Lambda=\fld{K}\Gamma/I$ such that $\Ld\al$ is uniserial for
every arrow $\al$ in $\Gamma_1$. Suppose $p$ is a mast for $U$ and
let $\alpha_1$ be the first arrow of $p$. Let $\mathcal{A}=\left\{
\Ld\gamma p\;|\; \gamma\in \Gamma_1 \right\}$. Any two members of
$\mathcal{A}$ are comparable; i.e.,  for $\gamma_1, \gamma_2 \in
\Gamma_1$, either $\Ld\gamma_1 p\subseteq \Ld\gamma_2 p$ or
$\Ld\gamma_2 p\subseteq \Ld\gamma_1 p$, since $\Ld \alpha_1$ is
uniserial. Hence there exists a greatest element in $\mathcal{A}$,
say $\Ld \gamma p$. Notice that $\Ld \gamma p$ can be zero. This
happens when $Jp=0$.\\
Case 1: There is no arrow leaving $e:=\s(p)$ except $\al_1$. Here
$\Ld e$ is uniserial, we have $U=\Lambda e/{\Lambda \gamma p}$ and
since $U$ is not projective, $\gamma p\neq 0$. Therefore $\al
(U)\leq 2\leq m$ by Proposition~\ref{indec}. \\
Case 2: There are arrows $\bt_1,\dotsc, \bt_l, \delta_{l+1},
\dotsc, \delta_n$ leaving $e$ except $\alpha_1$. Assume $(\bt_j,
e)\dt p$ ($1\leq j \leq l$) and $\delta_t e$ ($l+1\leq t\leq n$) are
non-routes on $p$. Note that $n < m$, since $\Ld$ is
$m$-multiserial. Let $b_j=\bt_j - \sum_{i\in I(\bt_j, e)}
k_i(\bt_j, e) v_i(\bt_j, e)$ and $b_t= \delta_t$.  If $\Ld\gamma
p=0$,  then $U=\Lambda e/{\sum_{i=1}^{n} \Lambda b_i}$. Otherwise
$U=\Lambda e /{(\sum_{i=1}^{n} \Lambda b_i +\Lambda\gamma p)}$.
Let
  $0\ra DTrU \xrightarrow{\:\:f\:\:} \bigsqcup_{i\in
I}B_i\xrightarrow{\:\:g\:\:} U\ra 0$ be an almost split sequence.
By Proposition~\ref{irrembeding}, all the induced irreducible maps
$g_i\colon B_i\ra U$ are epimorphisms. By
Proposition~\ref{mono}(iii), $\al (U)\leq \len \soc DTrU$. But
by~\cite[Proposition IV.1.11]{ars}, we know that $\soc DTrU \cong
P_1/JP_1$ where $P_1\ra \Lambda e \ra U\ra 0$ is a minimal
projective presentation of $U$. Therefore $\al (U)\leq m$.
\end{proof}

The following proposition gives more precise information in certain
situations; it follows from the proof of the above theorem.

\begin{prop}
Suppose $U$ is a non-projective uniserial module with mast $p$
over a left $m$-multiserial algebra $\Ld$.  Then
\begin{itemize}
\item[(i)]If $m=1$ (i.e.~if $\Ld$ is left serial) then $\al(U)\le 2$.
\item[(ii)]If there is only one arrow leaving $\s(p)$, then
$\al(U)\leq 2$.
\item[(iii)] If $m=2$ (for example, if $\Ld$ is a
left biserial algebra), and $Jp=0$, then $\al(U)=1$.
\end{itemize}
\end{prop}
\begin{proof}
The parts (i) and (ii) follow directly from the proof of the above
 theorem. As to part (iii),
let $\al_1$ be the first arrow of $p$. Then there is an
arrow $\bt\neq \al_1$ starting at $e=\s(p)$ (otherwise, following the
 above proof again, $U=\Ld e$ would
 be projective). Thus, either $(\bt,
e)\dt p$ or $\bt$ is a non-route on $p$. If $(\bt, e)\dt p$, then
$U=\Ld e/{\Ld b}$, where $b:=\bt -\sum_{i\in I(\bt, e)} k_i(\bt,
e) v_i(\bt, e)$. If $\bt$ is a non-route, then $U=\Ld e/{\Ld b}$,
where $b:=\bt$. In both cases then, $\al (U) =1$
by~\cite[Proposition~V.6.3]{ars}, because the image of $\Ld
f\xrightarrow{\:\: .b \:\:}\Ld e$ is not in $J^2e$, where
$f=t(\bt)$.
\end{proof}

\end{document}